\documentclass{amsart}
\usepackage{graphics}
\usepackage{amsmath,amsthm,amscd,amssymb}
\usepackage{epsfig}
\usepackage[all,cmtip]{xy}

\newcommand{\Q}{{\mathcal Q}}
\newcommand{\N}{\mathcal{N}}
\renewcommand{\mod}{\mathrm{mod\;}}

\newcommand{\gop}{{\bar{g}}}
\newcommand{\fop}{{\bar{f}}}
\newcommand{\eop}{{\bar{\varepsilon}}}
\newcommand{\alphaop}{{\bar{\alpha}}}
\def\gg{{\mathfrak u}}
\def\ff{{\mathfrak p}}
\def\ggop{{\bar{\mathfrak u}}}
\def\ffop{{\bar{\mathfrak p}}}
\def\xx{{\mathfrak x}}

\newcommand{\A}{\mathbf{\rm A}}

\newcommand{\e}{\varepsilon}

\newcommand{\Ext}{\operatorname{Ext}}
\newcommand{\Hom}{\operatorname{Hom}}
\newcommand{\Ker}{\operatorname{Ker}}
\newcommand{\Ima}{\operatorname{Im}}
\newcommand{\HH}{\operatorname{HH}}
\newcommand{\gr}{{\operatorname{gr}}}
\newcommand{\op}{{\operatorname{op}}}
\newcommand{\kar}{\operatorname{char}}
\newcommand{\MaxSpec}{\operatorname{MaxSpec}}
\newcommand{\rad}{\operatorname{rad}}
\newcommand{\soc}{\operatorname{soc}}

\newtheorem{Theorem}{Theorem}[section]
\newtheorem{Lemma}[Theorem]{Lemma}
\newtheorem{Proposition}[Theorem]{Proposition}
\newtheorem{Corollary}[Theorem]{Corollary}
\newtheorem{Definition}[Theorem]{Definition}
\theoremstyle{definition}
\newtheorem*{remark}{Remark}

\begin{document}

\topmargin .7cm \oddsidemargin 1.5cm \evensidemargin 1.5cm

\title[Hochschild cohomology and support varieties \dots]{Hochschild cohomology and support varieties for tame Hecke algebras}

\thanks{The first author acknowledges support through a Leverhulme Early
Career Fellowship}
\author{Sibylle Schroll}
\address{Department of Mathematics, University of
Leicester, University Road, Leicester LE1 7RH, UK }
\email{schroll@mcs.le.ac.uk}
\author{Nicole Snashall}
\address{Department of Mathematics, University of
Leicester, University Road, Leicester LE1 7RH, UK }
\email{N.Snashall@mcs.le.ac.uk}

\begin{abstract}

We give a basis for the Hochschild cohomology ring of tame Hecke
algebras. We then show that the Hochschild cohomology ring modulo nilpotence
is a finitely generated algebra of Krull dimension 2, and
describe the support varieties of modules for these algebras.

\end{abstract}

\date{\today}

\maketitle

\section*{Introduction}

\parskip7pt
\parindent0pt

Hecke algebras play an important role in representation theory. They
arise as deformations of the group algebras of finite Coxeter groups
and appear as endomorphism algebras of induced representations of
finite or $p$-adic Chevalley groups. They also give rise to the
Kazhdan-Lusztig polynomials which appear in the expression of the
canonical basis in terms of the natural basis of the Hecke algebra.
In particular, the Hecke algebras of type $A$ (which arise as a
deformation of the group algebras of the symmetric group) have been
well studied. A complete classification of the representation type
of the blocks of the Hecke algebras of type A was obtained in
\cite{Erdmann-Nakano}. For the tame  Hecke algebras of type $A$, it
was shown in \cite{Jost} that there are precisely two Morita
equivalence classes of blocks and they are represented by $H_q(S_4)$
and the principal block of $H_q(S_5)$ with $q = -1$. Moreover, they
are in the same derived equivalence class. This follows from
\cite{Membrillo-Hernandez} and see also
\cite{Bocian-Holm-Skowronski}, since they are generalized Brauer
tree algebras of the same type.

The algebra $H_q(S_5)$ with $q = -1$ is a Koszul symmetric special
biserial algebra.  Special biserial algebras occur in many aspects
of representation theory and they are necessarily of finite or tame
representation type. Symmetric special biserial algebras occur, for
example, as Hopf algebras associated to infinitesimal groups whose
principal block is tame \cite{Farnsteiner-Skowronski1,
Farnsteiner-Skowronski2}, in the representation theory of
$U_q(sl_2)$ \cite{Patra, Suter, Xiao}, of Drinfeld doubles of
generalized Taft algebras \cite{Erdmann-Green-Snashall-Taillefer}
(and see also \cite{Snashall-Taillefer-1}) and as socle deformations
of the latter \cite{Snashall-Taillefer-2}. For the symmetric special
biserial algebras in \cite{Snashall-Taillefer-1} and
\cite{Snashall-Taillefer-2}, the Hochschild cohomology ring modulo
nilpotence has been shown to be a finitely generated algebra of
Krull dimension 2. In this paper we show that the same phenomenon
holds for the tame Hecke algebras, so the question naturally arises
as to whether or not the Hochschild cohomology ring modulo
nilpotence of any symmetric special biserial algebra is a finitely
generated algebra of Krull dimension at most 2.


In this paper we let $K$ be an algebraically closed field and study
the Hochschild cohomology ring of the algebra $A$, which, in the
case where the characteristic of $K$ is not 2, is precisely the
principal block of $H_q(S_5)$ with $q = -1$. The algebra $A$ is a
symmetric special biserial algebra of tame representation type.
Furthermore, it is Koszul with radical cube zero, and so is one of
the algebras recently classified in \cite{Benson}; it corresponds to
the algebra associated to the simply laced extended Dynkin diagram
of type $\tilde{T}_2$ (we adopt the notation commonly used in
articles in physics journals; note that in \cite{Benson} and
\cite{Erdmann-Solberg} this diagram is denoted by $\tilde{Z}_1$). It
is well-known that the Hecke algebras of finite type arise from
Brauer tree algebras, and thus \cite{Erdmann-Holm} shows that they
are periodic algebras and describes their Hochschild cohomology
ring. Hence the structure of the Hochschild cohomology ring for
Hecke algebras of both tame and finite type is now known.
Specifically, in Theorem~\ref{Krull}, we prove that the Hochschild
cohomology ring modulo nilpotence of a Hecke algebra has Krull
dimension 1 if the algebra is of finite type and has Krull dimension
2 if the algebra is of tame type. For wild Hecke algebras, there has
so far been very little progress in the study of Hochschild
cohomology or support varieties.

The Hochschild cohomology ring modulo nilpotence was used in
\cite{Snashall-Solberg} to define the concept of a support variety for
any finitely generated module over a finite-dimensional algebra. It
was shown in \cite{Erdmann-Holloway-Snashall-Solberg-Taillefer}
that, under certain reasonable finiteness conditions
on the algebra, these support varieties have many of the analogous
properties to those satisfied by modules over a group algebra for a finite
group or over a cocommutative Hopf algebra. In addition, when these
conditions hold, the Hochschild cohomology ring is itself a finitely
generated algebra, and so the finiteness conjecture of \cite{Snashall-Solberg}
holds concerning the Hochschild cohomology ring modulo nilpotence.
We show that these finiteness conditions hold for the
tame Hecke algebras, which has also been independently shown by
Erdmann and Solberg in \cite{Erdmann-Solberg}. We then consider the
consequences for support varieties of modules over $A$.

\smallskip

The paper is structured as follows. In
section~\ref{bimoduleResolution}, we give a minimal projective bimodule
resolution of $A$, and in section~\ref{basis of HH} a basis for each
of the Hochschild cohomology groups $\HH^n(A)$ for $n \geq 0$. This
extends the results of \cite{Erdmann-Schroll}, where the dimensions
of the Hochschild cohomology groups were calculated for the Hecke
algebra $H_q(S_4)$ with $q=-1$. The advantage of considering $A$
over the algebras in \cite{Erdmann-Schroll} is that $A$ is a Koszul
algebra, and so has a linear projective resolution as an $A$-$A$-bimodule.
In addition, in section~\ref{moregeneral}, we give
the minimal projective bimodule resolution of a more general family of
algebras which contains the algebra $A$.
In section~\ref{fg1andfg2}, we recall the finiteness conditions {\bf Fg1.} and {\bf Fg2.} of \cite{Erdmann-Holloway-Snashall-Solberg-Taillefer} and use them to give new results on the support varieties of finitely generated $A$-modules.
Finally, (in Theorem~\ref{Krull}), we show that the Hochschild cohomology ring of $A$ modulo nilpotence is a finitely generated commutative algebra of Krull dimension 2.

\smallskip

The principal block of $H_q(S_5)$ for
$q = -1$ was given by quiver and relations in \cite{Erdmann-Nakano}.
Let $K$ be an algebraically closed field (with no restrictions placed on
$\kar K$). Let $\Q$ be the quiver
$$\xymatrix{
1\ar@(ul,dl)[]_{\varepsilon}\ar@<1ex>[r]^{\alpha} &
2\ar@(ur,dr)[]^{\bar{\varepsilon}}\ar@<1ex>[l]^{\bar{\alpha}} \\ }$$
and let $I$ be the ideal of $K\Q$ generated by $\{ \alpha \eop, \e
\alpha, \alphaop \e, \eop \alphaop, \e^2 - \alpha \alphaop, \eop^2 -
\alphaop \alpha \}$. Let $A = K\Q/I$. Then $A$ is the principal block of
$H_q(S_5)$ with $q = -1$ in the case where $\kar K$ is not 2.
The Hochschild cohomology ring
of $A$ is given by $\HH^*(A) = \Ext^*_{A^e}(A,A) =
\oplus_{n=0}^\infty \Ext^n_{A^e}(A,A)$ with the Yoneda product, where
$A^e = A^{\op} \otimes_k A$ denotes the enveloping algebra of $A$.

We denote the trivial path at the vertex $i$ by $e_i$. We write
paths from left to write. For any arrow $a$ in the quiver $\Q$, we
write ${\bf o}(a)$ for the trivial path corresponding to the origin
of $a$ and ${\bf t}(a)$ for the trivial path corresponding to the
terminus of $a$. Thus ${\bf o}(\varepsilon) = e_1 = {\bf
o}(\alpha)$, ${\bf t}(\varepsilon) = e_1$ and ${\bf t}(\alpha) =
e_2$ {\it etc.} There is an algebra isomorphism $b : A \rightarrow
A$ induced by the involution given by $e_1 \mapsto e_2$, $e_2
\mapsto e_1$, $\alpha \mapsto \alphaop$, $\alphaop \mapsto
\alpha$, $\e \mapsto \eop$ and $\eop \mapsto \e$.

\section{A minimal projective bimodule resolution}\label{bimoduleResolution}

The terms of the minimal projective bimodule resolution $(R_\bullet,
\delta)$ of $A$ can be calculated following \cite{Happel}.
We denote by $P_{ij}$ the projective indecomposable $A$-$A$-bimodule
$Ae_i \otimes e_j A$, for $i,j = 1,2$.
The projectives $R_n$ are
periodic of period $4$. We have $R_0 = P_{11} \oplus P_{22}$
and, for $4k +i \geq 1$,
\begin{equation}\label{terms}
R_{4k+i} \cong
\left\{\begin{array}{ll} P_{11}^{2k+1} \oplus
P_{12}^{2k+i} \oplus P_{21}^{2k+i} \oplus P_{22}^{2k+1}  & \mbox{if
$i=0,1$}
 \cr &\cr
 P_{11}^{2k+i-1} \oplus P_{12}^{2k+2}
\oplus P_{21}^{2k+2} \oplus P_{11}^{2k+i-1}  & \mbox{if $i=2,3$.}
\end{array} \right.
\end{equation}

We follow the approach of \cite{Green-Hartman-Marcos-Solberg} in defining the
minimal projective $A$-$A$-bimodule resolution of the
Koszul algebra $A$. As the starting point, we recall that
\cite{Green-Solberg-Zacharia} explains how to recursively define sets -
which we will denote by ${\mathcal G}^n$ rather than by $f^n$ - for each
$n \geq 0$, in order to give an explicit construction of a minimal projective
resolution $(Q_{\bullet}, d)$ of the right $A$-module $A/J(A)$ where $J(A)$ denotes the
Jacobson radical of $A$. These sets have the following properties.
\begin{enumerate}
\item[(i)] For each $n \geq 0$, $Q_n = \oplus_{x \in {\mathcal G}^n}{\bf
t}(x)A$.
\item[(ii)] For each $x \in {\mathcal G}^n$, there are unique elements $r_j
\in K\Q$ such that $x = \sum h_j^{n-1}r_j$ where the sum is over all
elements $h_j^{n-1} \in {\mathcal G}^{n-1}$.
\item[(iii)] For each $n \geq 1$, using the decomposition of (ii), for $x \in
{\mathcal G}^n$, the map $d_n : Q_n \to Q_{n-1}$ is given by ${\bf t}(x)\lambda
\mapsto \sum_jr_j{\bf t}(x)\lambda$.
\end{enumerate}
The minimal projective bimodule resolution of a Koszul algebra which was
given in \cite{Green-Hartman-Marcos-Solberg} uses these same sets of
\cite{Green-Solberg-Zacharia} in its construction.

Therefore,
for our algebra $A$, we start by defining the set ${\mathcal G}^0$ to be the set of vertices of
$\Q$, labelled so that
${\mathcal G}^0=\{ g^0_1 = e_1, \gop^0_1 =e_2 \}$. Then, for $n \geq 1$, we will recursively define
sets ${\mathcal G}^n = \{g_i^n, \gop_i^n, f_j^n ,\fop_j^n \}_{i,j}$ in $K\Q$
for appropriate indices $i$ and $j$ and such that $|{\mathcal G}^n| =
2(n+1)$ for all $n \geq 0$, so that these sets contain the requisite information to define
the minimal projective $A$-$A$-bimodule resolution of $A$.
The elements of the sets ${\mathcal G}^n$ will be uniform elements of $K\Q$, that is, for each $x \in
{\mathcal G}^n$ there are vertices $v, w$ such that $x = vxw$. We
set ${\bf o}(x) = v$ and ${\bf t}(x) = w$.

\begin{Definition}
For $n \geq 1$, we define elements:
{\small{
\begin{equation*}
\begin{array}{lllll}
g_1^n &=&
\left \{ \begin{array}{lll}
g_1^{n-1} \e - f_1^{n-1} \alphaop \;\;\;\;\;\;\;\;\;\;\; & \mbox{if $n=4k$ or $n=4k+2$} \\
g^0_1\e & \mbox{if $n=1$}\\
g_1^{n-1} \e - f_2^{n-1} \alphaop \;\;\;\;\;\;\;\;\;\;\; &
\mbox{if $n=4k+1$ with $k \geq 1$, or $n=4k+3$}\\
\end{array}\right. \\ \\
g_{2l}^n &=&
     \left\{ \begin{array}{lll} g_{2l}^{n-1} \e - f_{2l+1}^{n-1}
     \alphaop & \;\;\;\;\; & \mbox{$1 \leq l <k$ if $n=4k$
     or $1 \leq l \leq k$ if $n=4k+2$}\\
     g_{2l}^{n-1} \e + f_{2l-1}^{n-1}
     \alphaop & \;\;\;\;\; & \mbox{$1 \leq l \leq k$ if $n=4k+1$
     or $n=4k+3$}\\
     \end{array}\right. \\ \\
g_{2l+1}^n &=&
      \left\{ \begin{array}{lll} g_{2l+1}^{n-1} \e + f_{2l}^{n-1} \alphaop &
       \;\;\;\;\; & \mbox{$1 \leq l <k$ if $n=4k$ or $1 \leq l \leq k$ if
      $n=4k+2$} \\
       g_{2l+1}^{n-1} \e - f_{2l+2}^{n-1} \alphaop &
       \;\;\;\;\; & \mbox{$1 \leq l <k$ if $n=4k+1$ or $1 \leq l \leq k$ if
      $n=4k+3$} \\
      \end{array}\right. \\ \\
g_{2k}^{4k} &=& g_{2k}^{4k-1} \e \\ \\
g_{2k+1}^{n} &=&
            \left\{ \begin{array}{lll}  f_{2k}^{n-1} \alphaop  &
       \;\;\;\;\; & \mbox{ if $n=4k$}\\
            g_{2k+1}^{n-1} \e &
       \;\;\;\;\; & \mbox{ if $n=4k+1$ with $k \geq 1$}\\
            \end{array}\right. \\ \\
g_{2k+2}^{4k+3} &=& f_{2k+1}^{4k+2} \alphaop \\ \\
\end{array} $$
$$ \begin{array}{lllll}
f_1^n &=&
\left\{\begin{array}{ll}
  g_1^0\alpha & \mbox{if $n=1$}\\
  f_1^{n-1} \eop + g_1^{n-1} \alpha \;\;\;\;\;\;\;\;\;\;\;
  & \mbox{if $n=4k+1$ with $k \geq 1$, or $n=4k+3$} \\
  \end{array}\right.\\ \\
f_{2l+1}^n &=&
      \left\{ \begin{array}{lll} f_{2l+1}^{n-1} \eop - g_{2l+2}^{n-1}
      \alpha & \;\;\;\;\; & \mbox{$0
      \leq l <k$ if $n=4k$ or
      $n=4k+2$} \\
      f_{2l+1}^{n-1} \eop + g_{2l}^{n-1}
      \alpha & \;\;\;\;\; & \mbox{$1
      \leq l <k$ if $n=4k+1$ or $1 \leq l \leq k$ if
      $n=4k+3$} \\
      \end{array}\right. \\ \\
f_{2l+2}^n &=&
      \left\{ \begin{array}{lll}
      f_{2l+2}^{n-1} \eop + g_{2l+1}^{n-1} \alpha & \;\;\;\;\; & \mbox{$0
      \leq l <k$ if $n=4k$ or $n=4k+2$} \\
      f_{2l+2}^{n-1} \eop - g_{2l+3}^{n-1} \alpha & \;\;\;\;\; &
      \mbox{$0 \leq l < k$ if $n=4k+1$ or $n=4k+3$} \\
      \end{array}\right. \\ \\
f_{2k+1}^{n} &=&
         \left\{ \begin{array}{lll}
         g_{2k}^{n-1} \alpha  &
       \;\;\;\;\; & \mbox{ if $n=4k+1$ with $k \geq 1$} \\
        f_{2k+1}^{n-1} \eop  &
       \;\;\;\;\; & \mbox{ if $n=4k+2$}\\
        \end{array}\right. \\ \\
f_{2k+2}^{n} &=&
         \left\{ \begin{array}{lll}
        g_{2k+1}^{n-1} \alpha  &
       \;\;\;\;\; & \mbox{ if $n=4k+2$}\\
         f_{2k+2}^{n-1} \eop  &
       \;\;\;\;\; & \mbox{ if $n=4k+3$}.\\
        \end{array}\right. \\ \\
\end{array}
\end{equation*}
}}
Applying the bar involution we define $\gop_i^n = b(g^n_i)$ and
$\fop_j^n = b(f^n_j)$ and we obtain in this way all the elements of ${\mathcal G}^n$.
\end{Definition}

\begin{remark}
\begin{enumerate}
\item Note that $g^1_1 = \e, f^1_1 = \alpha, \gop^1_1 =
\eop$ and $\fop^1_1 = \alphaop$ so that ${\mathcal G}^1$ is the set
of arrows of $\Q$. Also, $g^2_1 = \e^2-\alpha\alphaop, f^2_1 =
\alpha\eop$ and $f^2_2 = \e\alpha$, and so $\gop^2_1 =
\eop^2-\alphaop\alpha, \fop^2_1 = \alphaop\e$ and $\fop^2_2 =
\eop\alphaop$. Thus ${\mathcal G}^2$ is a minimal set of uniform
elements which generate the ideal $I$.
\item It follows from above that
$$ \begin{array}{ll}
 {\mathcal G}^n = \{g_i^n, \gop_i^n,
f_j^n ,\fop_j^n \mid i=1, \ldots\ 2k+1, j =1, \ldots\, n-2k \} &
\mbox{for $n = 4k$ or $n = 4k+2$} \\
 {\mathcal G}^n = \{g_i^n, \gop_i^n,
f_j^n ,\fop_j^n \mid i=1, \ldots\ \frac{n+1}{2}, j =1, \ldots\,
\frac{n+1}{2} \} &
\mbox{for $n = 4k+1$ or $n = 4k+3$}
\end{array}$$
where $g_i^n \in e_1 K\Q e_1$, $f_j^n \in e_1 K\Q e_2$, $\fop_j^n \in
e_2 K\Q e_1$, and $\gop_i^n \in e_2 K\Q e_2$.
\item For each $n \geq 1$, the set ${\mathcal G}^n$ contains four
monomials of length $n$. These are the four subpaths of length $n$
of the path $(\e \alpha \eop \alphaop)^N$ for $N>>0$.
\end{enumerate}
\end{remark}

Thus $$ \begin{array}{lll}
R_n &=& \bigoplus_{x \in {\mathcal G}^n} A{\bf o}(x) \otimes{\bf t}(x) A\\ \\
&=& (\bigoplus_{i} A{\bf o}(g_i^n) \otimes {\bf t}(g_i^n) A) \oplus
(\bigoplus_{j} A{\bf o}(f_j^n) \otimes {\bf t}(f_j^n) A)\oplus
\\ \\
& & (\bigoplus_{i} A{\bf o}(\gop_i^n) \otimes {\bf t}(\gop_i^n) A) \oplus
(\bigoplus_{j} A{\bf o}(\fop_j^n) \otimes {\bf t}(\fop_j^n) A)
\end{array}$$
where, for each $i$, the summand $A{\bf o}(g_i^n) \otimes {\bf
t}(g_i^n) A = P_{11}$ and the summand $A{\bf o}(\gop_i^n) \otimes
{\bf t}(\gop_i^n) A = P_{22}$, and, for each $j$, the summand $A{\bf
o}(f_j^n) \otimes {\bf t}(f_j^n) A = P_{12}$ and the summand $A{\bf
o}(\fop_j^n) \otimes {\bf t}(\fop_j^n) A = P_{21}$.

\smallskip

Following \cite{Green-Hartman-Marcos-Solberg}, in order to define
the differential $\delta$ we need the following lemma so that we
have two different ways of expressing the elements of the set
${\mathcal G}^n$ in terms of the elements of the set ${\mathcal
G}^{n-1}$. Namely, we write the elements of ${\mathcal G}^n$ in the
form $\sum h_i^{n-1}p_i$ and in the form $\sum q_ih_i^{n-1}$ for
appropriate $p_i, q_i$ in the ideal of $K\Q$ which is generated by
all the arrows, and where the sums are over all $h_i^{n-1} \in
{\mathcal G}^{n-1}$. The proof of Lemma~\ref{Lemma:gn} is
straightforward and is omitted.

\begin{Lemma}\label{Lemma:gn}
For $n\geq 1$,
{\small{
$$
\begin{array}{lllll}
g_1^n &=&
\left\{\begin{array}{ll}
\e g^0_1 & \mbox{if $n=1$}\\
\e g_1^{n-1} - \alpha \fop_1^{n-1}
\;\;\;\;\;\;\;\;\;\;\; & \mbox{if $n \geq 2$} \\
\end{array}\right.\\ \\
g_{2l}^n &=& \e g_{2l+1}^{n-1} + \alpha \fop_{2l}^{n-1}
     \begin{array}{lll} \;\;\;\;\; & \mbox{where $1 \leq l <k$ if $n=4k$ or $1 \leq l \leq k$ if $n=4k+2$}\\
     \;\;\;\;\; & \mbox{or $1 \leq l \leq k$ if $n=4k+1$ or $n=4k+3$}\\
     \end{array}\\ \\
g_{2l+1}^n &=& \e g_{2l}^{n-1} - \alpha \fop_{2l+1}^{n-1}
      \begin{array}{lllll} \;\;\;\;\; & \mbox{where $1 \leq l <k$ if $n=4k$ or $1 \leq l \leq k$ if $n=4k+2$} \\
      \;\;\;\;\; & \mbox{or $1 \leq l <k$ if $n=4k+1$ or $1 \leq l \leq k$ if
      $n=4k+3$} \\
      \end{array} \\ \\
g_{2k}^{4k} &=&  \alpha \fop_{2k}^{n-1} \\ \\
g_{2k+1}^{n} &=& \e g_{2k}^{n-1}
       \;\;\;\;\; \mbox{ if $n=4k$, or if $n=4k+1$ with $k \geq 1$}\\ \\
g_{2k+2}^{4k+3} &=&  \alpha \fop_{2k+2}^{n-1}\\ \\
\end{array} $$
$$ \begin{array}{lllll}
f_1^n &=&  \left\{\begin{array}{ll}
\alpha \gop_1^0 & \mbox{if $n=1$}\\
\e f_2^{n-1} + \alpha \gop_1^{n-1} \;\;\;\;\;\;\;\;\;\;\;
& \mbox{if $n=4k+1$ with $k \geq 1$, or $n=4k+3$}\\
\end{array}\right.\\ \\
f_{2l+1}^n &=& \e f_{2l+2}^{n-1} + \alpha \gop_{2l+1}^{n-1}
      \begin{array}{lllll} \;\;\; & \mbox{where $0 \leq l <k$ if $n=4k$ or if $n=4k+2$} \\
      \;\;\; & \mbox{or $1 \leq l <k$ if $n=4k+1$ or $1 \leq l \leq k$ if $n=4k+3$} \\
      \end{array}\\ \\
f_{2l+2}^n &=& \e f_{2l+1}^{n-1} - \alpha \gop_{2l+2}^{n-1}
      \begin{array}{lllll} \;\;\;\;\; & \mbox{where $0
      \leq l <k$ if $n=4k$ or $n=4k+2$} \\
      \;\;\;\;\; & \mbox{or $0 \leq l < k$ if $n=4k+1$ or $n=4k+3$} \\
      \end{array} \\ \\
f_{2k+1}^{n} &=& \alpha \gop_{2k+1}^{n-1}
       \;\;\;\;\; \mbox{ if $n=4k+1$ with $k \geq 1$, or if $n=4k+2$} \\ \\
f_{2k+2}^{n} &=& \e f_{2k+1}^{n-1} \;\;\;\;\; \mbox{ if $n=4k+2$ or $n=4k+3$.} \\ \\
\end{array}
$$
 }}
and the analogous set of equations holds if we apply $b$ to the
equations above.
\end{Lemma}

We now introduce some additional notation to enable us to distinguish the summands of $R_n$. First, let
$$ \begin{array}{lllll}
\gg & := &  e_1 \otimes e_1 &\in& P_{11}= Ae_1 \otimes e_1A \\
\ggop & := &  e_2 \otimes e_2 &\in& P_{22}=Ae_2 \otimes e_2A \\
\ff & := &  e_1 \otimes e_2 &\in& P_{12}=Ae_1 \otimes e_2 A\\
\ffop & := &  e_2 \otimes e_1 &\in& P_{21}=Ae_2 \otimes e_1A. \\
\end{array}
$$
We add the superscript $n$ to indicate that each of these elements lies in the projective module $R_n$, and then add subscripts in order to
distinguish, within $R_n$, the different summands of the form $A e_i
\otimes e_jA$. Specifically, $\gg_i^n =   e_1 \otimes e_1$ lies in
the $i$th copy of $P_{11}$ as a summand of $R_n$, that is, in the
copy of $P_{11}$ corresponding to the element $g_i^n$ of ${\mathcal
G}^n$. Similarly, $\ff_j^n =   e_1 \otimes e_2$ lies in the $j$th
copy of $P_{12}$ as a summand of $R_n$, that is, in the copy of
$P_{12}$ corresponding to the element $f_j^n$ of ${\mathcal G}^n$.

The algebra isomorphism $b : A \to A$ extends to an algebra isomorphism
$b:A^e \to A^e$ where $b(a_1 \otimes a_2) = b(a_1) \otimes b(a_2)$ for $a_1 \in
A^{\op}, a_2 \in A$. Thus $b(\gg) = \ggop$ and $b(\ff) = \ffop$.

For better readability we omit in the following definition the
superscripts on all expressions except the left most one with the
understanding that the terms on the right hand side are in
$R_{n-1}$.

\begin{Definition}\label{differential}
We define an $A$-$A$-bimodule homomorphism $\delta_n: R_n \rightarrow
R_{n-1}$, for all $n\geq 1$, by
{\small{
$$
\begin{array}{lllll}
\gg_1^n &\mapsto &   \left\{ \begin{array}{lll} \gg_1 \e - \ff_1
\alphaop + \e \gg_1 - \alpha \ffop_1 & \;\;\;\;\; & \mbox{if $n=4k$
or $n=4k+2$} \\
\gg_1\e - \e\gg_1 & & \mbox{if $n=1$}\\
\gg_1 \e - \ff_2\alphaop - (\e \gg_1 - \alpha \ffop_1) & \;\;\;\;\; & \mbox{if
$n=4k+1$ with $k \geq 1$, or $n=4k+3$} \\
 \end{array}\right. \\ \\
\gg_{2l}^n &\mapsto &
     \left\{ \begin{array}{lllll} \gg_{2l} \e - \ff_{2l+1}
     \alphaop +  \e \gg_{2l+1} + \alpha \ffop_{2l}& \;\;\;\;\; & \mbox{$1 \leq l <k$ if $n=4k$
     or $1 \leq l \leq k$ if $n=4k+2$}\\
     \gg_{2l} \e + \ff_
     {2l-1}
     \alphaop - (\e \gg_
     {2l+1} + \alpha \ffop_{2l})& \;\;\;\;\; & \mbox{$1 \leq l \leq k$ if $n=4k+1$
     or $n=4k+3$}\\
     \end{array}\right. \\ \\
\gg_{2l+1}^n &\mapsto&
      \left\{ \begin{array}{lllll} \gg_{2l+1} \e + \ff_{2l} \alphaop +  \e \gg_{2l} - \alpha \ffop_{2l+1}&
       \;\;\;\;\; & \mbox{$1 \leq l <k$ if $n=4k$ or $1 \leq l \leq k$ if
      $n=4k+2$} \\
       \gg_{2l+1} \e - \ff_{2l+2} \alphaop - (\e \gg_{2l} - \alpha \ffop_{2l+1})&
       \;\;\;\;\; & \mbox{$1 \leq l <k$ if $n=4k+1$ or $1 \leq l \leq k$ if
      $n=4k+3$} \\
      \end{array}\right. \\ \\
\gg_{2k}^{4k} &\mapsto& \gg_{2k} \e + \alpha \ffop_{2k} \\ \\
\gg_{2k+1}^{n} &\mapsto&
            \left\{ \begin{array}{lllll}  \ff_{2k} \alphaop +  \e \gg_{2k}&
       \;\;\;\;\; & \mbox{ if $n=4k$}\\
            \gg_{2k+1} \e - \e \gg_{2k}&
       \;\;\;\;\; & \mbox{ if $n=4k+1$ with $k \geq 1$}\\
            \end{array}\right. \\ \\
\gg_{2k+2}^{4k+3} &\mapsto& \ff_{2k+1} \alphaop - \alpha \ffop_{2k+2}\\ \\
\end{array}$$
$$\begin{array}{lllll}
\ff_1^n &\mapsto& \left\{\begin{array}{ll}
   \gg_ 1 \alpha - \alpha \ggop_1 & \mbox{if $n=1$}\\
   \ff_1 \eop + \gg_ 1 \alpha - \e \ff_1 - \alpha \ggop_1 \;\;\;\;\;\;\;\;\;\;\;
   & \mbox{if $n=4k+1$ with $k \geq 1$, or $n=4k+3$} \\
   \end{array}\right.\\ \\
\ff_{2l+1}^n &\mapsto&
      \left\{ \begin{array}{lllll} \ff_{2l+1} \eop - \gg_{2l+2}
      \alpha  +  \e \ff_{2l+2} + \alpha \ggop_{2l+1}& \;\;\; & \mbox{$0
      \leq l <k$ if $n=4k$ or if
      $n=4k+2$} \\
      \ff_{2l+1} \eop + \gg_{2l}
      \alpha - (\e \ff_{2l+2} + \alpha \ggop_{2l+1})& \;\;\; & \mbox{$1
      \leq l <k$ if $n=4k+1$ or $1 \leq l \leq k$ if
      $n=4k+3$} \\
      \end{array}\right. \\ \\
\ff_{2l+2}^n &\mapsto&
      \left\{ \begin{array}{lllll}
      \ff_{2l+2} \eop + \gg_{2l+1} \alpha + \e \ff_{2l+1} - \alpha \ggop_{2l+2} & \;\;\;\;\; & \mbox{$0
      \leq l <k$ if $n=4k$ or $n=4k+2$} \\
      \ff_{2l+2} \eop - \gg_{2l+3} \alpha - (\e \ff_{2l+1} + \alpha \ggop_{2l+2})& \;\;\;\;\; &
      \mbox{$0 \leq l < k$ if $n=4k+1$ or $n=4k+3$} \\
      \end{array}\right. \\ \\
\ff_{2k+1}^{n} &\mapsto&
         \left\{ \begin{array}{lllll}
         \gg_{2k} \alpha  - \alpha \ggop_{2k+1}&
       \;\;\;\;\; & \mbox{ if $n=4k+1$ with $k \geq 1$} \\
          \ff_{2k+1} \eop  +  \alpha \ggop_{2k+1}&
       \;\;\;\;\; & \mbox{ if $n=4k+2$}\\
        \end{array}\right. \\ \\
\ff_{2k+2}^{n} & \mapsto&
         \left\{ \begin{array}{lllll}
        \gg_{2k+1} \alpha +  \e \ff_{2k+1} &
       \;\;\;\;\; & \mbox{ if $n=4k+2$}\\
         \ff_{2k+2} \eop - \e \ff_{2k+1} &
       \;\;\;\;\; & \mbox{ if $n=4k+3$}\\
        \end{array}\right. \\ \\
\end{array}
$$
}} together with the correspondences given by applying $b$, that is, for any
of the above correspondences $\xx \mapsto
\delta(\xx)$ we also have $b(\xx) \mapsto b(\delta(\xx))$.
\end{Definition}

The next result is now immediate from
\cite[Theorem 2.1]{Green-Hartman-Marcos-Solberg}.

\begin{Proposition}
The complex $(R_\bullet, \delta)$ is a minimal projective $A$-$A$-bimodule resolution of $A$, where $\delta_0 : R_0 \to A$ is the multiplication map.
\end{Proposition}

\subsection{A more general class of special biserial algebras}\label{moregeneral}

The algebra $A$ belongs to the wider class of special biserial algebras
whose quiver is $\Q$ but where we replace the relations by $\alpha \eop = \e
\alpha = \alphaop \e = \eop \alphaop = 0$, $\e^r - (\alpha
\alphaop)^s = 0$ and $\eop^r - (\alphaop \alpha)^s=0$, where $r, s$
are integers with $r \geq 2, s \geq 1$. Let $I'$ be the ideal of $K\Q$
generated by these relations and denote the resulting algebra by
$\Lambda = K\Q / I'$. We note that $\Lambda$ is Koszul if and only if $r=2$
and $s=1$, whence $\Lambda = A$. Although $\Lambda$ is, in general, not a
Koszul algebra, we may still give an explicit minimal
$\Lambda$-$\Lambda$-projective bimodule resolution of $\Lambda$.

The terms of the minimal projective bimodule resolution $(\tilde{R}_\bullet, \partial)$
of $\Lambda$ follow the same pattern as those of the corresponding resolution of $A$ given in~(\ref{terms}), that is, for all $r\geq 2$ and $s\geq 1$, we have
$\tilde{R}_0 = \tilde{P}_{11} \oplus \tilde{P}_{22}$ and, for $4k +i \geq 1$,
\begin{equation*}
\tilde{R}_{4k+i} \cong
\left\{\begin{array}{ll} \tilde{P}_{11}^{2k+1} \oplus
\tilde{P}_{12}^{2k+i} \oplus \tilde{P}_{21}^{2k+i} \oplus \tilde{P}_{22}^{2k+1}  & \mbox{if
$i=0,1$}
 \cr &\cr
 \tilde{P}_{11}^{2k+i-1} \oplus \tilde{P}_{12}^{2k+2}
\oplus \tilde{P}_{21}^{2k+2} \oplus \tilde{P}_{11}^{2k+i-1}  & \mbox{if $i=2,3$}
\end{array} \right.
\end{equation*}
where $\tilde{P}_{ij} = \Lambda e_i \otimes e_j\Lambda$ for $i, j  = 1, 2$.
By abuse of notation, we write $\gg = e_1 \otimes e_1 \in \tilde{P}_{11}$,
$\ggop = e_2 \otimes e_2 \in \tilde{P}_{22}$,
$\ff = e_1 \otimes e_2 \in \tilde{P}_{12}$ and $\ffop = e_2 \otimes e_1 \in \tilde{P}_{21}$, and denote by $b$ the algebra homomorphism $\Lambda^e \to \Lambda^e$ induced by the involution $e_1 \leftrightarrow e_2, \varepsilon \leftrightarrow \bar{\varepsilon}, \alpha \leftrightarrow \bar{\alpha}$ on $\Lambda$. We keep the same conventions on the use of
superscripts and subscripts as those preceding Definition~\ref{differential}, and make the following definition of
$\partial_n$; again we omit the superscripts on the right hand side
with the understanding that all the terms are in $\tilde{R}_{n-1}$.

\begin{Definition}
We define a $\Lambda$-$\Lambda$-bimodule homomorphism $\partial_n: \tilde{R}_n \rightarrow
\tilde{R}_{n-1}$, for all $n\geq 1$, by
{\small {
$$
\begin{array}{lllll}
\gg_1^n &\mapsto & \left\{
\begin{array}{lllll}
 \gg_1 \e^{r-1} - \ff_1 \alphaop (\alpha
\alphaop)^{s-1} + \e^{r-1} \gg_1 - (\alpha \alphaop)^{s-1} \alpha
\ffop_1 & \;\;\;\;\; & \mbox{if $n=4k$
or $n=4k+2$} \\
\gg_1 \e^{r-1} - \e^{r-1} \gg_1 & \;\;\;\;\; & \mbox{if $n=1$} \\
 \gg_1 \e^{r-1} - \ff_2 \alphaop (\alpha
\alphaop)^{s-1} - (\e^{r-1} \gg_1 - (\alpha \alphaop)^{s-1} \alpha
\ffop_1) & \;\;\;\;\; & \mbox{if $n=4k+1$ with $k \geq 1$,
or $n=4k+3$} \\
 \end{array}\right. \\ \\
\gg_{2l}^n &\mapsto &
     \left\{ \begin{array}{lllll} \gg_{2l} \e - \ff_{2l+1}
     \alphaop  (\alpha
\alphaop)^{s-1} +  \e^{r-1} \gg_{2l+1} + \alpha \ffop_{2l}&
\;\;\;\;\; & \mbox{$1 \leq l <k$ if $n=4k$
   or $1 \leq l \leq k$ if $n=4k+2$}\\
     \gg_{2l} \e^{r-1} + \ff_{2l-1}
     \alphaop - (\e^{r-1} \gg_{2l+1} + \alpha \ffop_{2l})& \;\;\;\;\; & \mbox{$1 \leq l \leq k$ if
     $n=4k+1$
 or if $n=4k+3$}\\
     \end{array}\right. \\ \\
\gg_{2l+1}^n &\mapsto&
      \left\{ \begin{array}{lllll} \gg_{2l+1} \e^{r-1} + \ff_{2l} \alphaop (\alpha
\alphaop)^{s-1}+  \e \gg_{2l} -  (\alpha \alphaop)^{s-1} \alpha
\ffop_{2l+1}&
       \;\; & \hspace{-.5cm} \mbox{$1 \leq l <k$ if $n=4k$ or $1 \leq l \leq k$ if
      $n=4k+2$} \\
       \gg_{2l+1} \e - \ff_{2l+2} \alphaop (\alpha \alphaop)^{s-1}- (\e \gg_{2l} -
       (\alpha \alphaop)^{s-1} \alpha \ffop_{2l+1})&
       \;\; & \hspace{-.5cm} \mbox{$1 \leq l <k$ if $n=4k+1$ or $1 \leq l \leq k$ if
      $n=4k+3$} \\
      \end{array}\right. \\ \\
\gg_{2k}^{4k} &\mapsto& \gg_{2k} \e + \alpha \ffop_{2k} \\ \\
\gg_{2k+1}^{n} &\mapsto&
            \left\{ \begin{array}{lllll}  \ff_{2k} \alphaop +  \e \gg_{2k}&
       \;\;\;\;\; & \mbox{ if $n=4k$}\\
            \gg_{2k+1} \e - \e \gg_{2k}&
       \;\;\;\;\; & \mbox{ if $n=4k+1$ with $k \geq 1$}\\
            \end{array}\right. \\ \\

\gg_{2k+2}^{4k+3} &\mapsto& \ff_{2k+1} \alphaop - \alpha \ffop_{2k+2}\\ \\
\end{array} $$
$$ \begin{array}{lllll}
\ff_1^n &\mapsto&  \left\{ \begin{array}{lllll} \gg_1 \alpha -
\alpha
\ggop_1 & \;\;\;\;\; & \mbox{if $n=1$} \\
 \ff_1 \eop^{r-1} + \gg_1 \alpha - \e^{r-1} \ff_2 - \alpha \ggop_1
& \;\;\;\;\; & \mbox{if $n=4k+1$ with $k \geq 1$ or $n=4k+3$} \\
\end{array}\right.
\\ \\
\ff_{2l+1}^n &\mapsto&
      \left\{ \begin{array}{lllll} \ff_{2l+1} \eop - \gg_{2l+2}
      \alpha (\alpha \alphaop)^{s-1} +  \e^{r-1} \ff_{2l+2} + \alpha \ggop_{2l+1}& \;\;\; & \mbox{$0
      \leq l <k$ if $n=4k$ or if
      $n=4k+2$} \\
      \ff_{2l+1} \eop^{r-1} + \gg_{2l}
      \alpha - (\e^{r-1} \ff_{2l+2} + \alpha \ggop_{2l+1})& \;\;\; & \mbox{$1
      \leq l <k$ if $n=4k+1$ or $1 \leq l \leq k$ if
      $n=4k+3$} \\
      \end{array}\right. \\ \\
\ff_{2l+2}^n &\mapsto&
      \left\{ \begin{array}{lllll}
      \ff_{2l+2} \eop^{r-1} + \gg_{2l+1} \alpha + \e \ff_{2l+1} - (\alpha \alphaop)^{s-1} \alpha \ggop_{2l+2} & \;\;\;\;\; & \mbox{$0
      \leq l <k$ if $n=4k$ or if $n=4k+2$} \\
      \ff_{2l+2} \eop - \gg_{2l+3} \alpha (\alpha \alphaop)^{s-1} - (\e \ff_{2l+1} - (\alpha \alphaop)^{s-1} \alpha \ggop_{2l+2})& \;\;\;\;\; &
      \mbox{$0 \leq l < k$ if $n=4k+1$ or $n=4k+3$} \\
      \end{array}\right. \\ \\
\ff_{2k+1}^{n} &\mapsto&
         \left\{ \begin{array}{lllll}
         \gg_{2k} \alpha  - \alpha \ggop_{2k+1}&
       \;\;\;\;\; & \mbox{ if $n=4k+1$ with $k \geq 1$} \\
          \ff_{2k+1} \eop  +  \alpha \ggop_{2k+1}&
       \;\;\;\;\; & \mbox{ if $n=4k+2$}\\
        \end{array}\right. \\ \\
\ff_{2k+2}^{n} & \mapsto&
         \left\{ \begin{array}{lllll}
        \gg_{2k+1} \alpha +  \e \ff_{2k+1} &
       \;\;\;\;\; & \mbox{ if $n=4k+2$}\\
         \ff_{2k+2} \eop^{r-1} - \e \ff_{2k+1} &
       \;\;\;\;\; & \mbox{ if $n=4k+3$}\\
        \end{array}\right. \\ \\
\end{array}
$$
}}
together with the correspondences given by applying $b$, that is,
for any $\xx \mapsto \partial(\xx)$ given above, we also have $b(\xx) \mapsto b(\partial(\xx))$.
\end{Definition}

\begin{Proposition}
The complex $(\tilde{R}_\bullet, \partial)$ is a minimal projective
$\Lambda$-$\Lambda$-bimodule resolution of $\Lambda$, where $\partial_0 : \tilde{R}_0 \to \Lambda$ is the multiplication map.
\end{Proposition}

{\it Proof:} A straightforward verification shows that $(\tilde{R}_\bullet,
\partial)$ is a complex. The fact that this complex is exact follows from an analogous argument to
the one in \cite[Proposition 2.8]{Green-Snashall}. 
 $\Box$

\medskip

\section{A basis of the Hochschild cohomology ring}\label{basis of
HH}

In this section we return to the algebra $A$ and give a basis of the
Hochschild cohomology ring $\HH^*(A)$ in the case where $\kar K \neq 2$.
By~\cite{Erdmann-Schroll} the dimensions of the
Hochschild cohomology groups $\HH^n(A)$ are as follows
$$ \dim \HH^{4k+i}(A) = \left\{
\begin{array}{ll}
2k+3 & \mbox{$i =0$ and $k>0$ or $i =1,2$ and $k\geq 0 $} \\
2k+4 & \mbox{$i =3, k \geq 0$ } \\
\end{array}
\right. $$ and the set $\{e_1+e_2, \e, \eop,
\e^2, \eop^2\}$ is a $K$-basis of  $\HH^0(A)$.

We describe the elements of a basis of each Hochschild cohomology
group $\HH^n(A)$, for $n \geq 1$, in terms of cocycles in $\Hom
(R_n, A)$. We write $\delta_n$ for both the map $\delta_n \colon R_n \to
R_{n-1}$ and for the induced map $\Hom(R_{n-1}, A) \to \Hom(R_n, A)$.
Thus our basis is given in terms of a set of elements in
$\Ker \delta_{n+1}$ as a subset of $\Hom (R^n, A)$, such that the
corresponding cosets in $\Ker \delta_{n+1} / \Ima \delta_n$ form a
basis of $\HH^n(A) = \Ker \delta_{n+1} / \Ima \delta_n$. We keep the notation of
section~\ref{bimoduleResolution}. When describing a
cocycle in $\Hom(R_n, A)$, we simply write the images of the generators
$\gg^n_i, \bar{\gg}^n_i, \ff^n_j$ or $\bar{\ff}^n_j$ in $R_n$ where those
images are non-zero.

It is easy to verify that the following elements are indeed cocycles in
$\Hom(R_n, A)$. They are used to define a basis of the Hochschild cohomology
groups in Theorem~\ref{basis}.

\begin{enumerate}
\item[(i)] Suppose $i = 0$ and $k \geq 1$ or $i = 2$ and $k \geq 0$.
\begin{enumerate}
\item[(a)] Let $\phi^{4k+i} : \begin{cases}
\gg^{4k+i}_1 \mapsto e_1, \\
\bar{\gg}^{4k+i}_1 \mapsto e_2.
\end{cases}$
\item[(b)] For $1 \leq l \leq k$, let $\theta_l^{4k+i} :
\left \{ \begin{array}{ll}
\gg^{4k+i}_{2l} \mapsto e_1, \hspace*{1cm} & \gg^{4k+i}_{2l+1} \mapsto e_1,\\
\bar{\gg}^{4k+i}_{2l} \mapsto e_2, & \bar{\gg}^{4k+i}_{2l+1} \mapsto e_2.
\end{array}\right.$
\item[(c)] For $1 \leq l \leq k$, let $\psi _l^{4k+i}:
\gg^{4k+i}_{2l} \mapsto \varepsilon^2$.
\item[]
\end{enumerate}
\item[(ii)] Suppose $i= 1$ or $i=3$ and that $k \geq 0$.
\item[]
For $1 \leq l \leq 2k+1$ if $i=1$ and $1 \leq l \leq 2k+2$ if $i=3$,
let $$\chi_l^{4k+i}: \left \{ \begin{array}{ll} \gg^{4k+i}_{l}
\mapsto \varepsilon,\hspace*{1.5cm} &
\bar{\gg}^{4k+i}_{l} \mapsto \bar{\varepsilon},\\
\ff^{4k+i}_{l} \mapsto (-1)^{l+1}\alpha, & \bar{\ff}^{4k+i}_{l}
\mapsto (-1)^{l+1}\bar{\alpha}. \end{array}\right.$$\\
\end{enumerate}

\begin{Theorem}\label{basis}
Suppose $\kar K \neq 2$, and let $k \geq 0$ and $0 \leq i \leq 3$ with $4k+i > 0$.
\begin{enumerate}
\item For $i= 0$ or $2$, $\HH^{4k+i}(A)$ has basis $\{ \phi^{4k+i}, \varepsilon \phi^{4k+i}, \bar{\varepsilon} \phi^{4k+i}, \psi_l^{4k+i}, \theta_l^{4k+i} \mid 1 \leq l \leq k\}$.
\item For $i= 1$ or $3$, $\HH^{4k+i}(A)$ has basis $\{ \chi_l^{4k+i}, \varepsilon \chi_1^{4k+i}, \bar{\varepsilon} \chi_1^{4k+i} \mid 1 \leq l \leq m_i\}$, where $m_i = 2k+1$ if $i=1$ and $m_i = 2k+2$ if $i=3$.
\end{enumerate}
\end{Theorem}

We omit the proof since it is straightforward to check these maps represent
linearly independent non-zero elements in $\HH^{4k+i}(A)$.

\section{{\bf Fg1.} and {\bf Fg2.}}\label{fg1andfg2}

Throughout this section, suppose again that $\kar K$ is arbitrary.
In \cite{Snashall-Solberg}, the Hochschild cohomology ring was used to
define the support variety of a finitely generated module over any
finite-dimensional algebra. It
was shown in \cite{Erdmann-Holloway-Snashall-Solberg-Taillefer} that, when
certain (reasonable) finiteness conditions hold, then the support varieties
of \cite{Snashall-Solberg} share many of the analogous properties of support
varieties for finite group rings or co-commutative Hopf algebras. We start by
stating these finiteness conditions and showing that they hold in our setting.

Let $J(A)$ denote the Jacobson radical of $A$. Then
the Yoneda algebra or Ext algebra of $A$ is given by $E(A) = \Ext^*(A/J(A),
A/J(A))$ with the Yoneda product. We use the notation $E(A)^n = \Ext^n(A/J(A), A/J(A))$ for the $n$th-graded component of $E(A)$.
The finiteness conditions of
\cite{Erdmann-Holloway-Snashall-Solberg-Taillefer} are:
\begin{enumerate}
\item[{\bf Fg1.}]  There is a graded subalgebra $H$ of $\HH^*(A)$ such that
$H$ is a commutative Noetherian ring and $H^0 = \HH^0(A)$.
\item[{\bf Fg2.}] $E(A)$ is a finitely generated $H$-module.
\end{enumerate}

The algebra $A$ is Koszul, so we know from \cite{BGSS} that the
image of the natural ring homomorphism $\HH^*(\A) \to E(A)$ is
precisely the graded centre $Z_{\gr}(E(A))$ of $E(A)$, where
$Z_{\gr}(E(A))$ is the subring of $E(A)$ generated by all
homogeneous elements $z$ in $E(A)^n$ ($n \geq 0$) such that $xy =
(-1)^{nm} yx$ for all $y \in E(A)^m$. Moreover, $E(A)$ is given by
quiver and relations as $E(A) = K\Q / I^{\perp}$ where $I^{\perp} =
\langle\e^2+\alpha \alphaop, \eop^2 + \alphaop \alpha \rangle$ (from
\cite[Theorem 2.2]{GM}). Thus $E(A)$ is the preprojective algebra
associated to the simply laced extended Dynkin graph of type
$\tilde{T}_2$.

Now, it is easy to check that the elements
$$x= \e^2 +  \eop^2 $$
and
$$ z =  \e \alpha \eop \alphaop + \alpha \eop \alphaop \e + \eop
\alphaop \e \alpha  + \alphaop \e \alpha \eop $$ are both in
$Z_{\gr}(E(A))$.

\begin{Theorem}\label{EAfg}
As a left $K[x,z]$-module, $E(A)$ is finitely generated with
generating set
$$S = \{e_1, e_2, \e, \alpha, \eop, \alphaop, \e
\alpha, \eop \alphaop, \alpha \eop, \alphaop \e, \e \alpha \eop,
\eop \alphaop \e,
 \alphaop \e \alpha, \alpha \e \alphaop, \e \alpha \eop \alphaop, \eop \alphaop \e \alpha  \}.$$
\end{Theorem}

Before proving Theorem~\ref{EAfg}, we consider the conditions {\bf Fg1.} and {\bf Fg2.}.
The element $\phi^2 \in \HH^2(A)$ is a pre-image of $x$ and the element
$\theta_1^4 \in \HH^4(A)$ is a pre-image of $z$. Let $H$
be the graded subalgebra of $\HH^*(A)$ generated by $\HH^0(A), \phi^2$ and
$\theta_1^4$, so that $H$ is a pre-image of $K[x,z]$ in $\HH^*(A)$.
Then we have the following immediate consequence of Theorem~\ref{EAfg}.

\begin{Corollary}\label{cor:fg}
The conditions {\bf Fg1.} and {\bf Fg2.} hold for the algebra $A$ with
respect to the subring $H$ of $\HH^*(A)$.
\end{Corollary}

This result was independently shown by Erdmann and Solberg in
\cite{Erdmann-Solberg}, where they consider all radical cube zero weakly
symmetric algebras. We observe that Corollary~\ref{cor:fg} does not require
the explicit structure of either $Z_{\gr}(E(A))$ or $\HH^*(A)$, and that
it now follows from \cite{Solberg} that {\bf Fg1.} and
{\bf Fg2.} also hold with respect to $\HH^{even}(A) \supset H$.

In order to prove Theorem~\ref{EAfg}, we need the following lemma.

\begin{Lemma}\label{Lemmageneratingset}
For any  $h \in S$ and any arrow $x$ in $\Q$ we have $hx = \sum_{s
\in S} z_s s$ where $z_s \in Z_{\gr}(E(A))$.
\end{Lemma}

{\it Proof:} We will show the result for $x = \alpha$ and $x = \e$;
the cases $x = \alphaop$ and $x= \eop$ are analogous.

Clearly the result holds for $\alpha$, $\e \alpha$, $\eop \alphaop
\alpha, \alphaop \e \alpha$ and $ \eop \alphaop \e \alpha$ . We note
that $\alphaop \alpha = -xe_2$, $\e^2 \alpha = x \alpha$ and $\alpha
\eop \alphaop \alpha = - \e^2 \alpha \eop = - x \alpha \eop$.
Finally, $\e \alpha \eop \alphaop \alpha = - \e^3 \alpha \eop = - x
\e \alpha \eop$. Again the result immediately holds for $\e,
\alphaop \e$ and $\eop \alphaop \e$. Furthermore, $\e^2 = xe_1$,
$\alphaop \e^2 = x \eop \alphaop$, $\alpha \eop \alphaop \e = ze_1 -
\e \alpha \eop \alphaop$ and $\e \alpha \eop \alphaop  \e = z \e - x
\alpha \eop \alphaop$. $\Box$

\medskip

{\it Proof of Theorem~\ref{EAfg}:} A proof by induction will show that $S$
is a generating set of $E(A)$.
More precisely, we show by induction on $n$ that each element
of $E(A)^n$ can be expressed as a linear combination of elements in
$S$ with coefficients in $Z_{\gr}(E(A))$.

The base case holds since $e_1$ and $e_2$ are in $S$. Suppose that
the result holds for $E(A)^n$, that is, every $x \in E(A)^n$ is of
the form $x= \sum_{s \in S} z_s s$ with $z_s \in Z_{gr}(E(A))$.
Since $E(A)$ is Koszul, an element $x \in E(A)^{n+1}$ is a linear
combination of elements of the form $y \gamma$ where $y \in E(A)^n$
and $\gamma$ is an arrow of $\Q$. By induction $y =\sum_{s \in S} z_s
s$ with $z_s \in Z_{gr}(E(A))$ and thus $y \gamma = \sum_{s \in S}
z_s s \gamma$. By Lemma~\ref{Lemmageneratingset} we have $y =
\sum_{s \in S} z_s (\sum_{t \in S} z'_t t)$ with $z'_t \in
Z_{gr}(E(A))$ and therefore $y = \sum_{s \in S} \sum_{t \in S} (z_s
z'_t) t$. Thus $x$ has the required form. $\Box$

\medskip

Following \cite{Snashall-Solberg}, the variety of a finitely generated
$A$-module $M$ with respect to $H$ is given by
$$V_H(M) = \MaxSpec (H/A_H(M,M))$$
where $A_H(M,M)$ is the annihilator of $\Ext^*_A(M,M)$ as an $H$-module.
Moreover, $V_H(A/J(A)) = \MaxSpec H$.
Corollary~\ref{cor:fg} and
\cite[Theorems 2.5(b), 5.3]{Erdmann-Holloway-Snashall-Solberg-Taillefer} have
the following direct consequences for our algebra $A$.

\begin{Proposition}
\begin{enumerate}
\item $\HH^*(A)$ is a finitely generated $K$-algebra.
\item For $M \in \mod A$, the variety of $M$ is trivial if and only if $M$
is projective.
\item For $M \in \mod A$ and $M$ indecomposable, the variety of $M$ is a
line if and only if $M$ is periodic.
\end{enumerate}
\end{Proposition}

We now consider the simple $A$-modules $S_1$ and $S_2$ (corresponding
to the vertices 1 and 2 respectively of the quiver $\Q$). From
\cite[Proposition 3.4]{Snashall-Solberg},
we have that $V_H(A/J(A)) = V_H(S_1) \cup V_H(S_2)$ and $V_H(S_1) =
V_H(\rad (e_1A))$. Using the exact sequence $0 \to S_1 \to
\rad (e_1A) \to A/J(A) \to 0$, we have $V_H(A/J(A)) \subseteq V_H(S_1) \cup
V_H(\rad (e_1A)) = V_H(S_1)$. Thus
$V_H(A/J(A)) = V_H(S_1)$. A similar argument for $S_2$ proves the
following proposition.

\begin{Proposition}
$V_H(S_1) = V_H(S_2) = \MaxSpec H$.
\end{Proposition}

This is, in fact, a special case of the following result.

\begin{Proposition}
Suppose that $\tilde{A}$ is a selfinjective algebra and that $\tilde{H}$ is a graded subalgebra of $\HH^*(\tilde{A})$ containing $\HH^0(\tilde{A})$. Let $S$ be the simple module corresponding to the indecomposable projective $P$, so that $P/\rad P \cong S$. If $\soc P \cong S$ and if $\tilde{A}/J(\tilde{A})$ is a summand of $\rad P/\soc P$ then $V_{\tilde{H}}(S) = V_{\tilde{H}}(\tilde{A}/J(\tilde{A})) = \MaxSpec \tilde{H}$.
\end{Proposition}

Finally, we remark that it can be verified that
$Z_{\gr}(E(A))$ and $K[x,z]$ are isomorphic
as algebras. Since $A$ is Koszul, we know from
\cite{BGSS} that the Hochschild cohomology ring of $A$ modulo
nilpotence is isomorphic to $Z_\gr(E(A))$ modulo nilpotence. Hence
$\HH^*(A)/\N \cong K[x,z]$, where $\N$ denotes the ideal
of $\HH^*(A)$ generated by all nilpotent elements. In contrast, for
a Hecke algebra $B$ of finite type, the results of
\cite{Erdmann-Holm} show that the Hochschild cohomology ring
modulo nilpotence of $B$ is isomorphic to $K[y]$, where $y$ is in
degree $m \geq 1$ and $m$ is minimal such that $\Omega_{B^e}^m(B)
\cong B$ as $B$-$B$-bimodules. Thus we have the following
description of the Hochschild cohomology ring modulo nilpotence for
all Hecke algebras of tame and of finite type.

\begin{Theorem}~\label{Krull}
The Hochschild cohomology ring modulo nilpotence of a Hecke algebra
has Krull dimension 1 if the algebra is of finite type and has Krull
dimension 2 if the algebra is of tame type.
\end{Theorem}

\end{document}